\input amstex
\documentstyle{amsppt}
\hsize160mm \vsize220mm \loadbold

 \centerline{\bf THE DISTANCE RATIO METRIC ON SOME DOMAINS OF $\Bbb C$ PLANE}

\vskip 1cm

\centerline{ \bf Slavko Simic}
\bigskip
\centerline{ Mathematical Institute SANU, Kneza Mihaila 36, 11000
Belgrade, Serbia}
\bigskip
\centerline{E-mail: ssimic\@turing.mi.sanu.ac.rs}

 \vskip 1cm


\vskip 1cm

{\sevenrm{\bf Abstract} We prove a generalization for some $\Bbb
C$ domains of Gehring - Palka theorem on M\"obius transformations
regarding the distance ratio metric. Namely, we show that this
theorem is valid for arbitrary holomorphic mappings $f:\Bbb
H\to\Bbb H$ or $f:\Bbb D\to\Bbb D$ with the same Lipschitz
constant $2$.}
\bigskip
{\bf 1. Introduction}
\bigskip
For a subdomain $G \subset {\Bbb R}^n\,$ and for all $x,y\in G$
the distance ratio metric $j_G$ is defined as
$$
 j_G(x,y)=\log \left( 1+\frac{|x-y|}{\min \{d(x,\partial G),d(y, \partial G) \} } \right)\,,
$$
where $d(x,\partial G)$ denotes the Euclidean distance from $x$ to
$\partial G$. The distance ratio metric was introduced by F.W.
Gehring and B.P. Palka [3] and in the above simplified form by M.
Vuorinen [9]. As the "first approximation" of the quasihyperbolic
metric, it is frequently used in the study of hyperbolic type
metrics ([1],[2],[6],[10]) and geometric theory of functions.
\bigskip
For an open continuous mapping $f: G \to G'$ we consider the
following condition: there exists a constant $C \ge 1$ such that
for all $x,y\in G$ we have
$$
j_{G'}(f(x), f(y)) \le C j_G(x,y) \,,
$$
or, equivalently, that the mapping
$$  f:(G, j_G) \to (G',j_{G'}) $$
between metric spaces is Lipschitz continuous with the Lipschitz
constant $C\,.$
\bigskip
However, unlike other frequently used metrics (hyperbolic metric,
absolute ratio metric), the distance ratio metric $j_G$ is not
invariant under M\"obius transformations. Therefore, it is natural
to ask what the Lipschitz constants are for this metric under
conformal mappings
 or M\"obius transformations in higher dimension.
F.\,W. Gehring, B.\,P. Palka and B.\,G. Osgood proved that $j$
metric is not changed by more than a factor $2$ under M\"obius
transformations, see [2], [3]:

\bigskip
{\bf Theorem A} \ {\it If $G$ and $G'$ are proper subdomains of
${\Bbb R}^n$ and if $f$ is a M\"obius transformation of $G$ onto
$G'$, then for all $x,y\in G$
$$ j_{G'}(f(x),f(y))\leq 2j_G(x,y)$$.}

\vskip 1cm

\bigskip
It is an interesting problem to investigate Lipschitz continuity
of the distance-ratio metric under some other conformal mappings.

In this paper we show that analogous to the above results are
valid for {\it arbitrary holomorphic self-mappings} of a
half-plane or the open unit disk. Thereby Theorem A is
considerably generalized in these cases.
\bigskip

\vskip 1cm

{\bf 2. Results and proofs}
\bigskip
 We shall treat firstly the case of upper half-plane $\Bbb
H:=\{z| \Im z>0\}.$

A generalization of Theorem A for this domain is given by the
following
\bigskip
{\bf Theorem 1} \ {\it For any holomorphic mapping $f, \ f:\Bbb
H\to \Bbb H, $ and all $z,w\in \Bbb H$, we have
$$
j_{\Bbb H}(f(z),f(w))\le 2 j_{\Bbb H}(z,w). \eqno(1.1)
$$

The constant $2$ is best possible. }
\bigskip
{\bf Proof}
\bigskip
Denote $s:=\min\{\Im z, \Im w\}, \ S:=\min\{\Im f(z), \Im f(w)\},
\ T:=\max\{\Im f(z), \Im f(w)\}$.

Then
$$
j_{\Bbb H}(z,w)=\log\Bigl(1+{|z-w|\over s}\Bigr); \ \ j_{\Bbb
H}(f(z),f(w))=\log\Bigl(1+{|f(z)-f(w)|\over S}\Bigr).
$$
\bigskip

Main tool in the proof will be Julia's variant of the famous
Schwarz-Pick Lemma for the half-plane [cf 2].
\bigskip
{\bf Theorem B} \ {\it For any holomorphic mapping $f, \ f:\Bbb
H\to \Bbb H, $ and all $z,w\in \Bbb H$, we have}
$$
|{f(z)-f(w)\over f(z)-\overline {f(w)}}|\le |{z-w\over z-\overline
w}|.
$$

\bigskip

Applying this assertion, we get
$$
|{f(z)-\overline {f(w)}\over f(z)-f(w)}|^2-1\ge |{z-\overline
w\over z-w}|^2-1,
$$

that is

$$
|{f(z)-f(w)\over z-w}|^2\le{\Im f(z)\Im f(w)\over \Im z\Im w},
$$

since for any $x,y\in \Bbb C$, the identity
$$
|x-\overline y|^2- |x-y|^2=4\Im x\Im y
$$

holds.
\bigskip
Therefore,
$$
{|f(z)-f(w)|^2\over S^2}\le {|z-w|^2\over s^2}{T\over
S}={|z-w|^2\over s^2}(1+{T-S\over S}),
$$

i.e.,

$$
{|f(z)-f(w)|\over S}\le {|z-w|\over s}\sqrt{1+{|f(z)-f(w)|\over
S}}, \eqno (1.2)
$$

since $T-S=|\Im(f(z)-f(w))|\le |f(z)-f(w)|$.

\bigskip

Denoting ${|z-w|\over s}:=2X\in [0, \infty)$, the relation $(1.2)$
gives
$$
1+{|f(z)-f(w)|\over S}\le (X+\sqrt{1+X^2})^2.
$$

Hence,

$$
j_{\Bbb H}(f(z),f(w))=\log\Bigl(1+{|f(z)-f(w)|\over S}\Bigr)\le
2\log(X+\sqrt{1+X^2})\le 2\log (1+2X)=2j_{\Bbb H}(z,w),
$$

and the proof is done.

\vskip 1cm

To prove that the constant 2 is best possible, choose
$$
f_0(z)=a-{1\over b+z},
$$

where $a$ and $b$ are arbitrary real numbers.
\bigskip

Since $\Im f_0(z)={\Im z\over |b+z|^2}$, we conclude that $f_0$
maps the upper half-plane into itself.
\bigskip
A calculation of $j$ values along the line $\zeta\subset \Bbb H$
given by $\zeta:=\{ z=i-b+t, \ t\in \Bbb R\}$ with $w:=i-b$, gives

$$
j_\zeta(z, w)=\log(1+t); \ j_\zeta(f_0(z),f_0(w))=\log(1+{{t\over
\sqrt{1+t^2}}\over {1\over 1+t^2}})=\log(1+t\sqrt{1+t^2}).
$$
\bigskip
Hence,
$$
{j_\zeta(f_0(z),f_0(w))\over j_\zeta(z,
w)}={\log(1+t\sqrt{1+t^2})\over \log(1+t)},
$$

and this expression is tending to 2 as $t\to\infty$.
\bigskip
It follows that the constant 2 is best possible.

\vskip 2cm

A generalization of Theorem A for the unit disk case is given by
the following
\bigskip
{\bf Theorem 2} \ {\it For any holomorphic mapping $f, \ f:\Bbb
D\to \Bbb D, $ and all $z,w\in \Bbb D$, we have
$$
j_{\Bbb D}(f(z),f(w))\le 2 j_{\Bbb D}(z,w). \eqno(2.1)
$$
}
\bigskip
{\bf Proof}
\bigskip
Let $\max\{|z|, |w|\}=r$ and suppose that $|f(z)|\ge |f(w)|$. Then
$$
j_{\Bbb D}(z,w)=\log(1+\frac{|z-w|}{1-r}); \ j_{\Bbb
D}(f(z),f(w))=\log(1+\frac{|f(z)-f(w)|}{1-|f(z)|}).
$$
\bigskip
The proof is based on well-known Schwarz-Pick lemma, stated in the
following form
\bigskip
{\bf Theorem C} \ {\it Let $f$ be a holomorphic mapping of the
open unit disk $\Bbb D$ into itself. Then for any $z_1, z_2\in
\Bbb D$, we have
$$
|\frac{f(z_2)-f(z_1)}{1-\bar{f}(z_1)f(z_2)}|\le
|\frac{z_2-z_1}{1-\bar{z}_1z_2}|.
$$}
\bigskip

An application of this lemma gives
$$
|\frac{f(z)-f(w)}{1-\bar{f}(z)f(w)}|^{-2}-1\ge
|\frac{z-w}{1-\bar{z}w}|^{-2}-1,
$$
that is,
$$
|\frac{f(z)-f(w)}{z-w}|^2\le\frac{(1-|f(z)|^2)(1-|f(w)|^2)}{(1-|z|^2)(1-|w|^2)}
$$
$$
\le \frac{(1-|f(z)|^2)(1-|f(w)|^2)}{(1-r^2)^2},
$$

where we used the following identity for complex numbers $x,y$,
$$
|1-\bar{x}y|^2-|x-y|^2=(1-|x|^2)(1-|y|^2).
$$
\bigskip
Therefore,
$$
\frac{|f(z)-f(w)|}{1-|f(z)|}\le\frac{|z-w|}{1-r}\frac{\sqrt{(1+|f(z)|)(1+|f(w)|)}}{1+r}\sqrt{\frac{1-|f(w)|}{1-|f(z)|}},
$$
i.e.,
$$
\frac{|f(z)-f(w)|}{1-|f(z)|}\le\frac{|z-w|}{1-r}\frac{1+|f(z)|}{1+r}\sqrt{1+\frac{|f(z)-f(w)|}{1-|f(z)|}},
\eqno(2.2)
$$
since
$$
\frac{1-|f(w)|}{1-|f(z)|}=1+\frac{|f(z)|-|f(w)|}{1-|f(z)|}\le
1+\frac{|f(z)-f(w)|}{1-|f(z)|}.
$$
\bigskip
 To obtain an estimation for $|f(z)|, z\in \Bbb D$, suppose that $f(0)=a\in \Bbb D$.
Applying Theorem C with $z_1=0, z_2=z$, we get
$$
|z|\ge
\frac{|f(z)-a|}{|1-\bar{a}f(z)|}\ge\frac{|f(z)|-|a|}{1-|a||f(z)|},
$$
that is,
$$
|f(z)|\le\frac{|z|+|a|}{1+|a||z|}\le\frac{r+|a|}{1+|a|r}\eqno(2.3)
$$

By (2.3) we get $\frac{1+|f(z)|}{1+r}\le
\frac{1+|a|}{1+|a|r}:=2c(a,r)=2c$ and, denoting
$\frac{|z-w|}{1-r}:=X\in [0, +\infty)$, the inequality (2.2) gives
$$
1+\frac{|f(z)-f(w)|}{1-|f(z)|}\le (cX+\sqrt{1+c^2X^2})^2.
$$

Hence,
$$
\frac{j_\Bbb D(f(z),f(w))}{j_\Bbb
D(z,w)}\le\frac{2\log(cX+\sqrt{1+c^2X^2})}{\log(1+X)}.
$$
\bigskip
Now it should be notified that the function $g$,
$$
g(X):=cX+\sqrt{1+c^2X^2}-(1+X)
$$
is negative for $X\in (0,T)$, where
$T=\frac{2(1-c)}{2c-1}=\frac{2|a|r+1-|a|}{|a|(1-r)}$.
\bigskip
Since $X=\frac{|z-w|}{1-r}\le\frac{2r}{1-r}$, it follows that
$X\in (0,T)$, i.e., $\log(cX+\sqrt{1+c^2X^2})\le\log(1+X)$,
\bigskip
which provides the proof.

\vskip 1cm

{\bf Remark} \ We cannot claim in this case that the best possible
Lipschitz constant is $C^*= 2$ since we do not know an example of
analytic function $f_0, \ f_0:\Bbb D\to\Bbb D$ such that
$$
\sup_{z,w}\frac{j_\Bbb D(f_0(z),f_0(w))}{j_\Bbb D(z,w)}=2.
$$

This problem was recently treated in [11] where the following
estimation for $C^*=C^*(a)$ is proved
$$
1+|a|\le C^*\le\min\{2(1+|a|), \sqrt{5+2|a|+|a|^2}\}; \ a:=f(0).
$$
\bigskip
By Theorem 2 this is reduced to
$$
1+|a|\le C^*\le 2,
$$

but problem of determining best possible Lipschitz constant for
the unit disk remains open.

\vskip 2cm

 {\bf References}
\bigskip

[1] \ { G.D. Anderson, M.K. Vamanamurthy and M. Vuorinen}, {\it
    Conformal Invariants, Inequalities and Quasiconformal Maps}, John
    Wiley \& Sons, New York, 1997.
\bigskip

[2] \ {H.P. Boas},{\it Julius and Julia: Mastering the Art of the
Schwarz Lemma}, Amer. Math. Monthly 117, November 2010, 770-785.

\bigskip

[3] \ { F.\,W. Gehring and B.\,G. Osgood},{\it Uniform domains and
the quasihyperbolic metric}, J. Analyse Math. 36, 1979, 50-74.

\bigskip
[4] \ { F.\,W. Gehring and B.P. Palka},{\it Quasiconformally
homogeneous domains}, J. Analyse Math. 30, 1976, 172-199.
\bigskip

[5] \ {  P. H\"ast\"o, Z. Ibragimov, D. Minda, S. Ponnusamy and S.
K. Sahoo}, {\it Isometries of some hyperbolic-type path metrics,
and the hyperbolic medial axis}, In the tradition of Ahlfors-Bers,
IV, Contemporary Math. 432, 2007, 63--74.
\bigskip

[6] \ { M. Huang, S.Ponnusamy, H. Wang, and X. Wang,} {\it A
cosine inequality in the hyperbolic geometry}, Appl. Math. Lett.
23 (2010), no. 8, 887--891.

\bigskip
[7] \ { R. Kl\'en, M. Vuorinen, and X. Zhang}, {\it
Quasihyperbolic metric and M\"obius transformations}.-
 Proc. Amer. Math. Soc. (to appear)
 Manuscript 9pp, 
 {arXiv: 1108.2967 math.CV}.
 \bigskip

[8] \ { J. V\"ais\"al\"a,} {\it The free quasiworld. Freely
quasiconformal and related maps in Banach spaces. Quasiconformal
geometry and dynamics}, (Lublin, 1996), 55--118, Banach Center
Publ., 48, Polish Acad. Sci., Warsaw, 1999.

\bigskip

[9] \ { M. Vuorinen}, {\it Conformal invariants and quasiregular
mappings}, {J. Analyse Math.,} {\bf 45} (1985), 69--115.

\bigskip

[10] \ { M. Vuorinen,} {\it Conformal geometry and quasiregular
mappings}, Lecture Notes in Math. 1319, Springer-Verlag,
Berlin-Heidelberg, 1988.
\bigskip
[11] \ S. Simic, {\it Lipschitz continuity of the distance ratio
metric on the unit disk}, Filomat 27/8 (2013).

\end